\documentclass[11pt]{amsart}

\usepackage{amssymb}
\usepackage{amsmath}
\usepackage{mathrsfs}
\usepackage{url}
\usepackage{hyperref}
\input xy
\xyoption{all}
\usepackage{mathdots}
\usepackage{graphicx}  
\usepackage{tikz}
\usepackage{tkz-graph}

\newtheorem{thm}{Theorem}

\newtheorem{cor}[thm]{Corollary}
\newtheorem{prop}[thm]{Proposition}
\newtheorem{conj}[thm]{Conjecture}

\newtheorem{Definition}[thm]{Definition}
\newenvironment{definition}
  {\begin{Definition}\rm}{\end{Definition}}

\newtheorem{Example}[thm]{Example}
\newenvironment{example}
  {\begin{Example}\rm}{\end{Example}}
  
\newtheorem{Remark}[thm]{Remark}
\newenvironment{remark}
  {\begin{Remark}\rm}{\end{Remark}}

\title{Another proof of Wilmes' conjecture}

\author{Sam Hopkins}
\email{shopkins@mit.edu}
\address{Massachusetts Institute of Technology, Cambridge MA, 02139}

\subjclass[2010]{05E40, 05C25, 13D02, 06B30}
\keywords{Chip-firing, G-parking function, lcm-lattice}

\begin{document}

\begin{abstract}
We present a new proof of the monomial case of Wilmes' conjecture, which gives a formula for the coarsely-graded Betti numbers of the $G$-parking function ideal in terms of maximal parking functions of contractions of~$G$. Our proof is via poset topology and relies on a theorem of Gasharov, Peeva, and Welker~\cite{gasharov} that connects the Betti numbers of a monomial ideal to the topology of its $\mathrm{lcm}$-lattice.
\end{abstract}

\maketitle

\section{Introduction: the \texorpdfstring{$G$}{G}-parking function ideal} \label{sec:intro}

Wilmes' conjecture concerns the Betti numbers of a certain polynomial ideal closely related to the chip-firing game on a graph. Chip-firing on a graph has been studied in various contexts and under various names, including graphical parking functions~\cite{postnikov}, the Abelian sandpile model~\cite{perkinson}, and discrete Riemann-Roch theory~\cite{baker}. For a comprehensive introduction to the sandpile theory behind the problem, see~\cite{perkinson}, especially~\S7. We state the conjecture here concisely and without broader context; however, a few definitions are required first:

Let $G$ be an unoriented, connected graph with vertex set $V = \{v_1,\ldots,v_n\}$ and edge set $E$. Loops are not permitted, but multiple edges are allowed; formally, we represent the edges of $G$ by a finite set of labels $E$ together with a map $\varphi\colon E \to \binom{V}{2}$. For a subset $W \subseteq V$, let $G_W$ denote the subgraph induced on~$W$ by $G$. A \emph{connected partition}~$\pi = \{W_1,\ldots,W_k\}$ of $G$ is a partition of $V$ such that~$G_{W_i}$ is connected for all $i$. We use $|\pi| := k$ to denote the number of parts. Note that the set of connected partitions of $G$ is independent of the multiplicity of edges: it depends only on which vertices are adjacent. For a connected partition $\pi$ of $G$, let $G\vert_{\pi}$ be the contraction of $G$ to $\pi$; that is, let $G\vert_{\pi}$ be the graph obtained from $G$ by collapsing all the vertices in each part of $\pi$ into one vertex. A \emph{cut} of $G$ is a partition of $V$ into two parts. A cut is \emph{connected} if it is a connected partition.

From now on, we fix $v_n$ as the sink vertex of $G$. Let $R := \mathbf{k}[x_1,\ldots,x_{n-1}]$ be a polynomial ring in $n-1$ variables over a field $\mathbf{k}$. For a cut $C = \{U,W\}$ with~$v_n \in W$, define
\[ \mathbf{x}^{C} := \prod_{i=1}^{n-1} x_i^{d_U(v_i)},\]
where 
\[d_U(u) := \begin{cases}
|\{ e \in E\colon \varphi(e) = \{u,w\} \textrm{ with } w \in V\setminus U\}| & \textrm{if $u \in U$}, \\
0 &\textrm{otherwise}.\end{cases}\]
Define the monomial ideal~$I$ of $R$ by
\[ I := \langle \mathbf{x}^{C}\colon C \textrm{ is a connected cut of $G$}\rangle.\]
We call $I$ the \emph{$G$-parking function ideal}\footnote{The $G$-parking function ideal can be constructed in various other, equivalent ways. See~\cite{mania}, \cite{mohammadi}, \cite{manjunath}, or~\cite{dochtermann} for a proof that the $x^C$ are in fact the minimal generators of $I$.}, because a monomial $\mathbf{x} = \prod_{i=1}^{n-1} x_i^{c_i}$ is nonvanishing on $R/I$ exactly when $\sum_{i=1}^{n-1} c_i v_i$ is a $G$-parking function:

\begin{definition}
Set $\widetilde{V} := V \setminus \{v_n\}$. A \emph{$G$-parking function} $c = \sum_{i = 1}^{n-1} c_i v_i$ with respect to $v_n$ is an element of~$\mathbb{Z}\widetilde{V}$ such that for every non-empty subset~\mbox{$U \subseteq \widetilde{V}$}, there exists $v_i \in U$ with $0 \leq c_i < d_{U}(v_i)$. The parking functions inherit a partial order from $\mathbb{Z}\widetilde{V}$. A parking function $c'$ is \emph{maximal} if $c' \leq c$ for any parking function $c$ implies $c = c'$.
\end{definition}

The number of $G$-parking functions, and the number of maximal $G$-parking functions, are independent of the choice of sink vertex~\cite{benson}. We use~$\mathrm{mpf}(\Gamma)$ to denote the number of maximal parking functions of any connected graph $\Gamma$. The following conjecture about the coarsely graded Betti numbers\footnote{The Betti numbers are topological invariants of an $R$-module that can be read off from a minimal free resolution of that module. See~\cite{miller} for an introduction to Betti numbers.} of $R/I$ was stated in~\cite{wilmes} and again in~\cite{perkinson}:

\begin{conj} \label{con:wilmes} \emph{(Wilmes' conjecture)}
For all $i \geq 1$, we have
\[ \beta_i(R/I) = \sum_{|\pi| = i+1} \mathrm{mpf} (G\vert_{\pi}),\]
where the sum is over all connected partitions $\pi$ of $G$ with $i+1$ parts.
\end{conj}

Actually, Wilmes originally conjectured that this sum gives the Betti numbers of the \emph{topppling ideal}, a binomial ideal of which $I$ is a distinguished monomial initial ideal. The $G$-parking function ideal was introduced by Cori, Rossin, and Salvy~\cite{cori} and its resolutions have been studied by Postnikov and Shapiro~\cite{postnikov} and Manjunath and Sturmfels~\cite{manjunath2}. At any rate, Wilmes' conjecture has now been proven several times: by Mohammadi and Shokrieh~\cite{mohammadi} (for both the monomial and binomial cases), by Manjunath, Schreyer, and Wilmes~\cite{manjunath} (again for both the monomial and binomial cases), and by Dochtermann and Sanyal~\cite{dochtermann} (in the monomial case). Also, Mania~\cite{mania} has a proof for $\beta_{1}$ in the binomial case. We present a new proof here (in the monomial case) because we believe the ideas are significantly novel. The proof, via poset topology~\cite{wachs}, rests on a result by Gasharov, Peeva, and Welker~\cite{gasharov} that connects the Betti numbers of a monomial ideal to the topology of the $\mathrm{lcm}$-lattice of the ideal. The idea is to construct an ideal $J$, the \emph{connected cut-set ideal} of $G$, whose $\mathrm{lcm}$-lattice is dual to the connected partition lattice of $G$, and then to show that $I$ and~$J$ have the same Betti numbers. The dual connected partition lattice, $L^*_G$, is useful for our purposes because for any~$\pi \in L^*_G$, the M\"{o}bius function evaluated at $\pi$ is equal (up to sign) to the number of maximal parking functions on $G\vert_{\pi}$. We use one further ideal $K$, the \emph{oriented connected cut-set ideal} of $G$, to connect~$I$ and $J$. Example~\ref{ex:kite} illustrates all of the objects we attach to the graph $G$ in order to prove Conjecture~\ref{con:wilmes}: the dual connected partition lattice~$L^*_G$, its M\"{o}bius function, and the three ideals $I$, $J$, and $K$.

\begin{example} \label{ex:kite}
Let $G$ be the ``kite graph'' pictured below:
\begin{center}
\begin{tikzpicture}
\SetVertexMath
\GraphInit[vstyle=Art]
\SetUpVertex[MinSize=3pt]
\SetVertexLabel
\tikzset{VertexStyle/.style = {
shape = circle,
shading = ball,
ball color = black,
inner sep = 2pt
}}
\SetUpEdge[color=black]
\Vertex[LabelOut,Lpos=90,
Ldist=.1cm,x=0,y=1]{v_2}
\Vertex[LabelOut,Lpos=180,
Ldist=.1cm,x=-0.75,y=0]{v_1}
\Vertex[LabelOut,Lpos=0,
Ldist=.1cm,x=0.75,y=0]{v_3}
\Vertex[LabelOut,Lpos=270,
Ldist=.1cm,x=0,y=-1]{v_4}
\Edges(v_2,v_1,v_4,v_3,v_2)
\Edges(v_1,v_3)
\node at (0.8, -1) {$G$};
\node at(-0.5,0.7){$a$};
\node at(0,0.2){$b$};
\node at(-0.5,-0.7){$c$};
\node at(0.6,0.6){$d$};
\node at(0.6,-0.6){$e$};
\end{tikzpicture}
\end{center}
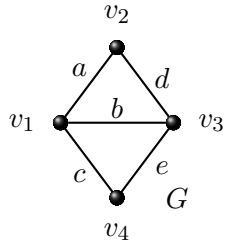
Figure~\ref{fig:lattice} shows the dual connected partition lattice $L^*_G$~(\S\ref{sec:bond}), whose atoms are labeled by generators of the $G$-parking function ideal~$I$~(\S\ref{sec:intro}), the connected cut-set ideal~$J$~(\S\ref{sec:cut}), and the oriented connected cut-set ideal~$K$~(\S\ref{sec:orientedcut}), and all of whose elements are labeled by their M\"{o}bius function values. (For $K$, we orient the edges  so that if $f \in E$ with $f_1 = v_i$ and $f_2 = v_j$, then $i < j$.) $\square$
\end{example}

\begin{figure} \label{fig:lattice}
\begin{tikzpicture}[scale=0.85]

\def \x {0}
\def \y {0}
{	
	\node at (\x,\y+1) [circle,draw,fill=black,inner sep=1pt]{};
	\node at (\x-0.75,\y) [circle,draw,fill=black,inner sep=1pt]{};
	\node at (\x+0.75,\y) [circle,draw,fill=black,inner sep=1pt]{};
	\node at (\x,\y-1) [circle,draw,fill=black,inner sep=1pt]{};
	\draw (\x,\y+1) -- (\x-0.75,\y) -- (\x,\y-1) -- (\x+0.75,\y) -- (\x,\y+1);
	\draw (\x-0.75,\y) -- (\x+0.75,\y);
}

\draw plot [smooth cycle] coordinates { (0,-1.25) (-1,0) (0,1.25) (1,0)};

\draw (-0.75,0.75) -- (-5,1.75);
\draw (-0.6,0.9) -- (-3,1.75);
\draw (-0.4,1.2) -- (-1,1.75);
\draw (0.4,1.2) -- (1,1.75);
\draw (0.6,0.9) -- (3,1.75);
\draw (0.75,0.75) -- (5,1.75);

\def \y {5}
\foreach \x in {-5,-3,-1,1,3,5}
{
	\node at (\x,\y+1) [circle,draw,fill=black,inner sep=1pt]{};
	\node at (\x-0.75,\y) [circle,draw,fill=black,inner sep=1pt]{};
	\node at (\x+0.75,\y) [circle,draw,fill=black,inner sep=1pt]{};
	\node at (\x,\y-1) [circle,draw,fill=black,inner sep=1pt]{};
	\draw (\x,\y+1) -- (\x-0.75,\y) -- (\x,\y-1) -- (\x+0.75,\y) -- (\x,\y+1);
	\draw (\x-0.75,\y) -- (\x+0.75,\y);
}

\draw (-5,6.5) -- (-4.5,7.5);
\draw (-4.5,6) -- (-3,7.5);

\draw (-3.25,6.25) -- (-4.25,7.5);
\draw (-2.5,6) -- (2,7.75);

\draw (-1.5,6) -- (-4,8);
\draw (-1,6.3) -- (-0.1,7.6);
\draw (-0.5,6) -- (4,7.75);

\draw (1,6.3) -- (0.1,7.6);
\draw (0.5,6) -- (-2,7.6);
\draw (1.5,6) -- (2,7.6);

\draw (3.25,6.25) -- (4.25,7.5);
\draw (2.5,6) -- (-2,7.75);

\draw (5,6.5) -- (4.5,7.5);
\draw (4.5,6) -- (3,7.5);

\node at (-5,2.75){\parbox[h]{1 in}{\begin{center}$x_1^3$\\$y^{a}y^{b}y^{c}$\\$z_1^{a}z_1^{b}z_1^{c}$\end{center}}};
\node at (-3,2.75){\parbox[h]{1 in}{\begin{center}$x_1^2x_2$\\$y^{b}y^{c}y^{d}$\\$z_1^{b}z_1^{c}z_1^{d}$\end{center}}};
\node at (-1,2.75){\parbox[h]{1 in}{\begin{center}$x_2^2$\\$y^{a}y^{d}$\\$z_2^{a}z_1^{d}$\end{center}}};
\node at (1,2.75){\parbox[h]{1 in}{\begin{center}$x_1x_3$\\$y^{c}y^{e}$\\$z_1^{c}z_1^{e}$\end{center}}};
\node at (3,2.75){\parbox[h]{1 in}{\begin{center}$x_2x_3^2$\\$y^{a}y^{b}y^{e}$\\$z_2^{a}z_2^{b}z_1^{e}$\end{center}}};
\node at (5,2.75){\parbox[h]{1 in}{\begin{center}$x_3^2$\\$y^{b}y^{d}y^{e}$\\$z_2^{b}z_2^{d}z_1^{e}$\end{center}}};

\draw plot [smooth cycle] coordinates {(-5,5-1.25) (-5+0.9,5) (-5,5+1.25) (-5.25,5+1) (-5.5+1,5)  (-5.25,5-1)};
\draw (-5.75,5) circle (0.25);

\draw plot [smooth cycle] coordinates { (-4,5) (-3,6.25) (-2.75,6) (-3.75,4.75)};
\draw plot [smooth cycle] coordinates { (-3.25,4) (-2.25,5.25) (-2,5) (-3,3.75)};

\draw plot [smooth cycle] coordinates { (-1.9,5.15) (-0.1,5.15) (-1,3.85)};
\draw (-1,6) circle (0.25);

\draw plot [smooth cycle] coordinates { (1.9,4.85) (0.1,4.85) (1,6.15)};
\draw (1,4) circle (0.25);

\draw plot [smooth cycle] coordinates { (4,5) (3,6.25) (2.75,6) (3.75,4.75)};
\draw plot [smooth cycle] coordinates { (3.25,4) (2.25,5.25) (2,5) (3,3.75)};

\draw plot [smooth cycle] coordinates {(5,5-1.25) (5-0.9,5) (5,5+1.25) (5.25,5+1) (5.5-1,5)  (5.25,5-1)};
\draw (5.75,5) circle (0.25);

\def \y {9}
\foreach \x in {-4.5,-2.25,0,2.25,4.5}
{
	\node at (\x,\y+1) [circle,draw,fill=black,inner sep=1pt]{};
	\node at (\x-0.75,\y) [circle,draw,fill=black,inner sep=1pt]{};
	\node at (\x+0.75,\y) [circle,draw,fill=black,inner sep=1pt]{};
	\node at (\x,\y-1) [circle,draw,fill=black,inner sep=1pt]{};
	\draw (\x,\y+1) -- (\x-0.75,\y) -- (\x,\y-1) -- (\x+0.75,\y) -- (\x,\y+1);
	\draw (\x-0.75,\y) -- (\x+0.75,\y);
}

\draw (-4.25,10.25) -- (-0.5,11.5);
\draw (-2,10.25) -- (-0.3,11.3);
\draw (0,10.3) -- (0,10.7);
\draw (4.25,10.25) -- (0.5,11.5);
\draw (2,10.25) -- (0.3,11.3);

\draw plot [smooth cycle] coordinates { (-4.75,8) (-3.75,9.25) (-3.5,9) (-4.5,7.75)};
\draw (-4.5-0.75,9) circle (0.25);
\draw (-4.5,10) circle (0.25);

\draw plot [smooth cycle] coordinates { (-1.25,9) (-2.25,10.25) (-2.5,10) (-1.5,8.75)};
\draw (-2.25-0.75,9) circle (0.25);
\draw (-2.25,8) circle (0.25);

\draw plot [smooth cycle] coordinates { (-0.9,9.15) (0.9,9.15) (0.9,8.85) (-0.9,8.85)};
\draw (0,10) circle (0.25);
\draw (0,8) circle (0.25);

\draw plot [smooth cycle] coordinates { (4.75,8) (3.75,9.25) (3.5,9) (4.5,7.75)};
\draw (4.5+0.75,9) circle (0.25);
\draw (4.5,10) circle (0.25);

\draw plot [smooth cycle] coordinates { (1.25,9) (2.25,10.25) (2.5,10) (1.5,8.75)};
\draw (2.25+0.75,9) circle (0.25);
\draw (2.25,8) circle (0.25);

\def \x {0}
\def \y {12}
{	
	\node at (\x,\y+1) [circle,draw,fill=black,inner sep=1pt]{};
	\node at (\x-0.75,\y) [circle,draw,fill=black,inner sep=1pt]{};
	\node at (\x+0.75,\y) [circle,draw,fill=black,inner sep=1pt]{};
	\node at (\x,\y-1) [circle,draw,fill=black,inner sep=1pt]{};
	\draw (\x,\y+1) -- (\x-0.75,\y) -- (\x,\y-1) -- (\x+0.75,\y) -- (\x,\y+1);
	\draw (\x-0.75,\y) -- (\x+0.75,\y);
}

\draw (0,13) circle (0.25);
\draw (0,11) circle (0.25);
\draw (0.75,12) circle (0.25);
\draw (-0.75,12) circle (0.25);

\node at (1,-1) [shape=rectangle,draw]{1};

\node at (-6,4) [shape=rectangle,draw]{-1};
\node at (-4,4) [shape=rectangle,draw]{-1};
\node at (-2,4) [shape=rectangle,draw]{-1};
\node at (2,4) [shape=rectangle,draw]{-1};
\node at (4,4) [shape=rectangle,draw]{-1};
\node at (6,4) [shape=rectangle,draw]{-1};

\node at (-5.25,8) [shape=rectangle,draw]{2};
\node at (-2.9,8) [shape=rectangle,draw]{2};
\node at (-0.8,8.25) [shape=rectangle,draw]{1};
\node at (2.9,8) [shape=rectangle,draw]{2};
\node at (5.25,8) [shape=rectangle,draw]{2};

\node at (1,13) [shape=rectangle,draw]{-4};

\end{tikzpicture}
\caption{The dual connected partition lattice $L^*_G$ for $G$ the ``kite graph'' of Example~\ref{ex:kite}. For each $\pi \in L^*_G$, the value of $\mu(\hat{0},\pi)$ appears in a square box beside $\pi$.}
\end{figure}
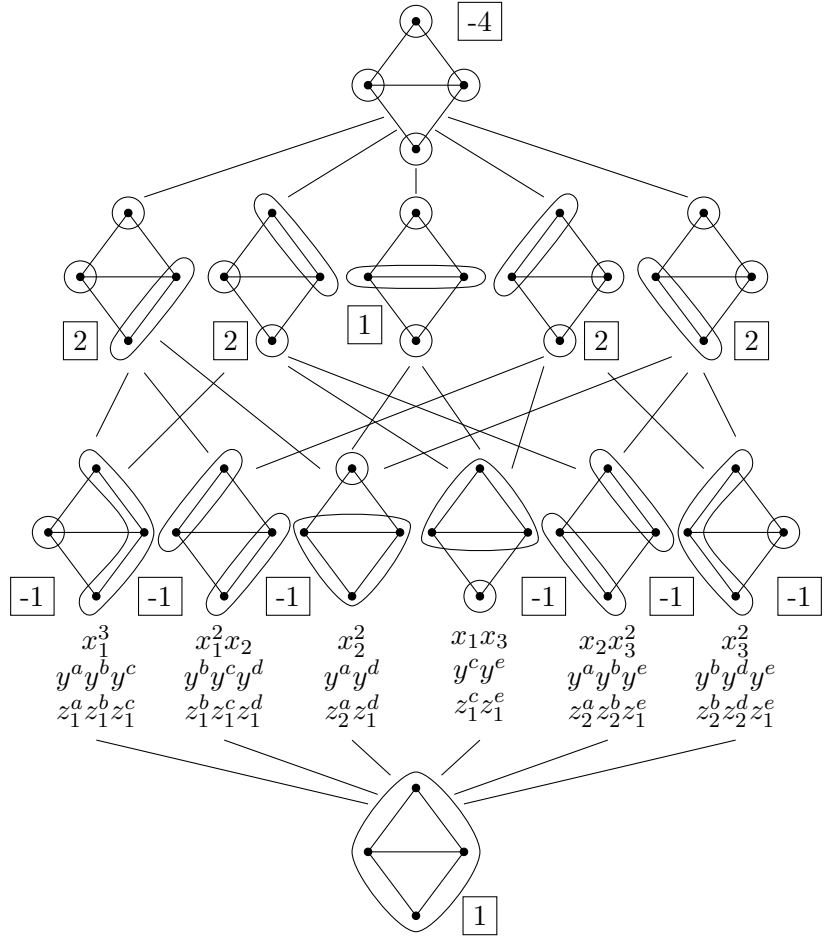

\begin{remark} Many of the ideas presented here are similar to those found in a preprint of Mohammadi and Shokrieh~\cite{mohammadi2} that was posted to the arXiv shortly after the first version of this note. Mohammadi and Shokrieh explain the matroidal structure behind the ideals considered in \S\ref{sec:cut} and \S\ref{sec:orientedcut}:  the ideal $J$ is the \emph{(unoriented) matroid ideal} of $G$ and $K$ is the \emph{graphic oriented matroid ideal}.\footnote{For a direct translation between this note and~\cite{mohammadi2}: our ideal $I$ is their ideal $\mathbf{M}_G^{v_n}$; our ideal $J$ is their ideal $\mathbf{Mat}_G$; and our ideal $K$ is their ideal $\mathbf{O}_G^{v_n}$.} The content of~\cite{mohammadi2} is much more comprehensive; for instance, it addresses also the binomial toppling ideal. Nevertheless, the emphasis on the $\mathrm{lcm}$-lattice here differentiates our approach from that of~\cite{mohammadi2}.
\end{remark}

\noindent {\bf Acknowledgements}: We thank David Perkinson for the helpful discussion, comments, and proofreading. We also thank the anonymous reviewer for help with the exposition and references.

\section{The connected partition lattice and its dual} \label{sec:bond}

\begin{definition}
The \emph{connected partition lattice} $L_G$ is the set of connected partitions of $G$ partially ordered by refinement: $\pi < \pi'$ if $\pi$ refines $\pi'$. We use $L^*_G$ to denote the lattice that is the order dual of $L_G$. We use $\hat{0}$ and $\hat{1}$ to denote the minimal and maximal elements of a lattice, respectively. So we have~$\hat{0} = \{ \{v_1\},\ldots,\{v_n\}\}$ and $\hat{1} = \{V\}$ for $L_G$, and vice-versa for $L^*_G$. 
\end{definition}

Note that $L_G$ is in fact the lattice of flats of the graphic matroid associated to $G$~\cite[\S5]{novik} and is thus geometric. The lattice $L^*_G$ is also geometric: dualizing preserves semi-modularity and that it is atomic follows from Proposition~\ref{prop:lcmdual} below. The M\"{o}bius function of the dual connected partition lattice can be computed in terms of maximal parking functions, as Proposition~\ref{prop:bondmob} explains. This proposition is~\cite[Proposition 5.3]{novik}.

\begin{prop} \label{prop:bondmob}
The M\"{o}bius function of $L^*_G$ is given by
 \[ \mu(\hat{0},\pi) = (-1)^{|\pi|-1} \mathrm{mpf} (G\vert_{\pi})\]
 for all $\pi \in L^*_G$.
 \end{prop}
 
\section{The connected cut-set ideal} \label{sec:cut}

\begin{definition}
For a cut~$\{U,W\}$, we define its \emph{cut-set} to be 
\[\{ e \in E\colon \varphi(e) = \{u,w\} \textrm{ with } u \in U, w \in W\}.\] 
A \emph{connected cut-set} of $G$ is the cut-set of a connected cut.
\end{definition}

Let $S :=  \mathbf{k}[y^{e}\colon e \in E]$ be a polynomial ring in~$|E|$ variables over $\mathbf{k}$. For~$F \subseteq E$, define
\[\mathbf{y}^{F} := \prod_{\substack{e \in F \\ }} y^{e}.\] 
Define the squarefree monomial ideal $J$ of~$R$ by
\[ J := \langle \mathbf{y}^F\colon F \subseteq E \textrm{ is a connected cut-set of $G$} \rangle.\]
We call $J$ the \emph{connected cut-set ideal} of $G$. It is clear that the generators above are minimal. Let $L_J$ be the $\mathrm{lcm}$-lattice of $J$: the lattice of least common multiples of minimal generators $\mathbf{y}^F$ ordered by divisibility ($\mathbf{y} \leq \mathbf{y}'$ if $\mathbf{y}$ divides $\mathbf{y}'$). 

 \begin{prop} \label{prop:lcmdual}
We have the isomorphism of lattices $L_J \simeq L^*_G$.
 \end{prop}
 
 \emph{Proof}: For a connected partition $\pi \in L^*_G$, define
 \[ F(\pi) := \{ e \in E\colon \varphi(e) = \{u,v\} \textrm{ with $u$ and $v$ not in the same part of $\pi$}\}.\]
 We claim that $\phi\colon \pi \mapsto \mathbf{y}^{F(\pi)}$ is a bijection between $L^*_G$ and $L_J$. An atom~$\pi$ of $L^*_G$ is just a connected cut, and in this case $F(\pi)$ is a connected cut-set, so clearly $\phi$ bijects between the atoms of $L^*_G$ and atoms of $L_J$. 
 
 Let $\pi,\pi' \in L^*_G$. Their join can be obtained as follows. First, find their common refinement as partitions:
 \[ \pi'' := \{ W \cap W'\colon W \in \pi, W' \in \pi',W \cap W' \neq \emptyset\}.\]
 But $\pi''$ in general is not a connected partition because $G_{W}$ for $W \in \pi''$ maybe be disconnected. Thus,
 \[ \pi \vee \pi' = \{ U\colon G_U \textrm{ is a connected component of some $G_{W}$ for $W \in \pi''$}\}.\]
 It is easy to see that $F(\pi \vee \pi') = F(\pi) \cup F(\pi')$: if $e \in F(\pi)$ or~$e \in F(\pi')$ with $\varphi(e) = \{u,w\}$, then $u$ and $v$ will not be in the same part of $\pi \vee \pi''$; on the other hand, if $e \notin F(\pi)$ and $e \notin F(\pi')$, then $u$ and $v$ will still be in the same part of $\pi \vee \pi''$. And then observe that $\mathbf{y}^{F(\pi)} \vee \mathbf{y}^{F(\pi')} = \mathbf{y}^{F(\pi) \cup F(\pi')}$. $\square$
 
Let us use $\mathrm{rk}(x)$ to denote the rank of an element $x$ of a lattice $L$. Since~$L_J$ is a geometric lattice, each interval $(\hat{0},\mathbf{y})$ in $L_J$ has the homotopy type of a wedge of spheres of dimension equal to~$\mathrm{rk}(\mathbf{y}) - 2$.\footnote{When we say that a poset has a topological property, we mean that the \emph{order complex} of that poset has that property. See~\cite{wachs} for the definition of the order complex of a poset.}  Thus, the following theorem of Gasharov, Peeva, and Welker~\cite[Theorem~2.1]{gasharov} lets us compute the Betti numbers of $S/J$ from the M\"{o}bius function of $L_J$. Actually these Betti numbers were computed and given a combinatorial interpretation already in~\cite{novik}; but the $\mathrm{lcm}$-lattice makes this computation straightforward.
 
 \begin{thm}
 The $i$th coarsely graded Betti numbers of $S/J$ are given by
 \[\beta_i(S/J) = \sum_{\substack{\mathbf{y} \in L_J \\ \mathbf{y} \neq \hat{0}}} \mathrm{dim} \; \widetilde{\mathrm{H}}_{i-2}( (\hat{0},\mathbf{y}); \mathbf{k})\]
 for all $i \geq 1$. Here $\widetilde{\mathrm{H}}_{i-2}( (\hat{0},\mathbf{y}); \mathbf{k})$ is the reduced homology of $(\hat{0},\mathbf{y})$.
 \end{thm}
 
 \begin{cor} \label{cor:cutbetti}
By the above theorem, the Betti numbers are
\begin{align*}
\beta_i(S/J) = \sum_{\substack{\mathbf{y} \in L_J \\ \mathrm{rk}(\mathbf{y}) = i}} |\mu(\hat{0},\mathbf{y})| = \sum_{\substack{\pi \in L_G \\ |\pi| = i+1}} \mathrm{mpf} (G\vert_{\pi})
\end{align*}
for all $i \geq 1$.
\end{cor}
\emph{Proof}: For any $\mathbf{y} \in L_J$, the interval $(\hat{0},\mathbf{y})$ has vanishing reduced homology in every dimension except dimension $\mathrm{rk}(\mathbf{y}) - 2$. Then by the Euler characteristic~\cite{wachs} we get~$\mathrm{dim} \; \widetilde{\mathrm{H}}_{\mathrm{rk}(\mathbf{y})-2}( (\hat{0},\mathbf{y}); \mathbf{k}) = |\mu(\hat{0},\mathbf{y})|$. The second equality follows from Propositions~\ref{prop:bondmob} and~\ref{prop:lcmdual}. $\square$

\section{The oriented connected cut-set ideal} \label{sec:orientedcut}

Let $T := \mathbf{k}[z_{1}^e, z_{2}^e\colon e \in E]$ be a polynomial ring in $2|E|$ variables over~$\mathbf{k}$. Recall that we have fixed $v_n$ as the sink of $G$. Let us also choose $e_1,e_2 \in V$ for each $e \in E$ so that $\varphi(e) = \{e_1,e_2\}$. In other words, let us fix an orientation of each edge of $G$. For a cut $C = \{U,W\}$ with~$v_n \in W$, define
\[\mathbf{z}^{C} := \prod_{ \substack{e \in E \\ e_1 \in U, e_2 \in W}} z_{1}^{e} \cdot \prod_{ \substack{e \in E \\ e_2 \in U, e_1 \in W}} z_{2}^{e} .\]
Define the squarefree monomial ideal~$K$ of $T$ by
\[ K := \langle \mathbf{z}^{C}\colon C \textrm{ is a connected cut of $G$}\rangle.\]
We call $K$ the \emph{oriented connected cut-set ideal} of $G$. Again, these generators are minimal. The oriented connected cut-set ideal serves as the bridge between the connected cut-set ideal and the $G$-parking function ideal, as Propositions~\ref{prop:orientedtopark} and~\ref{prop:orientedtocut} demonstrate.

\begin{prop} \label{prop:orientedtopark}
For each $1 \leq i \leq n$, fix an arbitrary $e^{i} \in E$, $\varepsilon^i \in \{1,2\}$ such that $e^{i}_{\varepsilon^{i}}=v_i$. Then the sequence
\[A := \{z_{\varepsilon^{i}}^{e^{i}} - z_{\delta}^{f}\colon 1 \leq i \leq n, f \in E, \delta \in \{1,2\}, f_\delta = v_i, \textrm{ and $f \neq e^i$}\} \cup \{z_{\varepsilon^n}^{e^n}\}\] 
is a permutable regular sequence on $T/K$. Further, $T / K \otimes_T T/ \langle A \rangle \simeq R/I$.
\end{prop}

\emph{Proof}: The isomorphism is clear from construction. We now prove that~$A$ is a permutable regular sequence. Our strategy is similar to the proof of the analogous fact for the polarization of a monomial ideal given in~\cite{swanson}. This proof depends crucially on a monotonicity property of the generators of~$K$ (so observe the connection with~\cite[\S5]{postnikov}).  Let~$\widetilde{A}$ be any subsequence of $A$. We can safely ignore $z_{\varepsilon^n}^{e^n}$ because it will never divide any generator of $K$. Thus it suffices to show that if~$z_{\varepsilon^{i}}^{e^{i}} - z_{\delta}^{f} \notin \widetilde{A}$ with~$z_{\varepsilon^{i}}^{e^{i}} - z_{\delta}^{f} \in A$, then~$z_{\varepsilon^{i}}^{e^{i}} - z_{\delta}^{f}$ is non-zero-divisor of~$T / K \otimes_T T/ \langle \widetilde{A} \rangle$. Without loss of generality suppose that~$\varepsilon^{i} = \delta = 1$; so in particular, $e^i_1 = f_1 = v_i$.

Let $\widetilde{K}$ be the ideal of $T/\langle \widetilde{A} \rangle$ obtained from the generators of $K$ by identifying all the $z_{\gamma}^{g}$ with $z_{\varepsilon^{j}}^{e^j}$ for $z_{\varepsilon^{j}}^{e^{j}} - z_{\gamma}^{g} \in \widetilde{A}$. Note $T / K \otimes_T T/ \langle \widetilde{A} \rangle=(T/\langle \widetilde{A} \rangle)/\widetilde{K}$. Write the ideal~$\widetilde{K}$ as
\[\widetilde{K} = z_1^{e^i}\widetilde{K}_1 + z_1^{f}\widetilde{K}_2 + z_1^{e^i}z_1^{f}\widetilde{K}_3 + \widetilde{K}_4,\] 
where the minimal monomial generators of~$\widetilde{K}_1$ do not involve $z_1^{f}$, the minimal monomial generators of~$\widetilde{K}_2$ do not involve~$z_1^{e^i}$, and the minimal monomial generators of~$\widetilde{K}_3$ and of $\widetilde{K}_4$ involve neither of these. 

Suppose $(z_1^{e^i} - z_1^{f})$ is a zero-divisor in $(T/\langle \widetilde{A} \rangle)/\widetilde{K}$. Then it lies in some associated prime ideal of $\widetilde{K}$, and hence $z_1^{e^i}$ and $ z_1^{f}$ also lie in this associated prime. So there exists some~$r \in (T/\langle \widetilde{A} \rangle)\setminus \widetilde{K}$ with $rz_1^{e^i} \in\widetilde{K}$ and $rz_1^{f} \in \widetilde{K}$.

This gives us~$r \in \widetilde{K}_1 + z_1^{f}\widetilde{K}_2 + z_1^{f}\widetilde{K}_3 + \widetilde{K}_4$ and~$r \in z_1^{e^i}\widetilde{K}_1 +\widetilde{K}_2 + z_1^{e^i}\widetilde{K}_3 + \widetilde{K}_4$. So we have,
\begin{align*}
r &\in (\widetilde{K}_1 + z_1^{f}\widetilde{K}_2 + z_1^{f}\widetilde{K}_3 + \widetilde{K}_4) \cap (z_1^{e^i}\widetilde{K}_1 +\widetilde{K}_2 + z_1^{e^i}\widetilde{K}_3 + \widetilde{K}_4) \\
&= \widetilde{K} + \widetilde{K}_1 \cap \widetilde{K}_2 + \widetilde{K_1} \cap z_1^{e^i} \widetilde{K}_3 + \widetilde{K}_2 \cap z_1^{f} \widetilde{K}_3.
\end{align*}
Minimal monomial generators of $\widetilde{K}_1 \cap \widetilde{K}_2$ are given by $\mathrm{lcm}(\mathbf{z}',\mathbf{z}'')$, where~$\mathbf{z}'$ ranges over minimal monomial generators of $\widetilde{K}_1$ and $\mathbf{z}''$ ranges over minimal monomial generators of $\widetilde{K}_2$. Let $\mathbf{z}'$ be a minimal generator of $\widetilde{K}_1$ and~$\mathbf{z}''$ be a minimal generator of~$\widetilde{K}_2$. Note that the monomial $z_1^{e^i}\mathbf{z}'$ corresponds to some cut $C_1 = \{U_1,W_1\}$ with $e^i_1,f_2 \in U_1$ and~$e^i_2, v_n \in W_1$, and similarly~$z_1^{f}\mathbf{z}''$ corresponds to a cut $C_2 = \{U_2, W_2\}$ with $e^i_1,e^i_2 \in U_2$ and~$f_2,v_n \in W_2$. Set~\mbox{$W := W_1 \cap W_2$} and~$U := U_1 \cup U_2$. Let $W_3 \subseteq W$ be such that $G_{W_3}$ is the connected component of $G_W$ that contains $v_n$. Set~\mbox{$C_3 := \{V\setminus W_3, W_3\}$}. Note that $C_3$ is connected: by definition, $G_{W_3}$ is connected, and $G_{U}$ is connected since $e^i_1$ and $e^i_2$ and $e^i_1$ and $f_2$ are adjacent. Finally, there must be an edge between a vertex in each connected component of~$G_W$ and a vertex in $G_{U}$ since $G$ as a whole is connected and these connected components are disconnected from one another. Next we claim that $\mathbf{z}^{C_3}$ (with the identifications of $z_{\gamma}^{g}$ and~$z_{\varepsilon^{j}}^{e^j}$ as above) divides $\mathrm{lcm}(\mathbf{z}',\mathbf{z}'')$: let $F_1, F_2, F_3$ be the cut-sets of~$C_1, C_2, C_3$; then we have~$F_3 \subseteq (F_1 \cup F_2)$ and $e^i,f \notin F_3$. So~$\widetilde{K}_1 \cap \widetilde{K}_2 \subseteq \widetilde{K}$. Similarly one can show $ \widetilde{K_1} \cap z_1^{e^i} \widetilde{K}_3$ and $\widetilde{K}_2 \cap z_1^{f} \widetilde{K}_3$ are subsets of $\widetilde{K}$, and therefore~$r \in \widetilde{K}$, a contradiction. Thus, $z_1^{e^i} - z_1^{f}$ cannot be a zero-divisor. $\square$

\begin{prop} \label{prop:orientedtocut}
The sequence $B := \{z_{1}^{e} - z_{2}^e \colon e\in E\}$ is a permutable regular sequence on $T/K$. Further, $T / K \otimes_T T/ \langle B \rangle \simeq S/J$.
\end{prop}

\emph{Proof}: This follows from~\cite[Corollary 2.7]{novik}. Alternatively, essentially the same proof as in the last proposition works again. $\square$

\begin{thm} For $i \geq 1$, $\displaystyle \beta_i(R/I) = \beta_i(T/K) = \beta_i(S/J) = \hspace{-0.2cm} \sum_{\substack{\pi \in L_G \\ |\pi| = i+1}} \hspace{-0.2cm} \mathrm{mpf} (G\vert_{\pi})$.
\end{thm}

\vspace{-0.3cm}

\emph{Proof}: Taking the quotient of $T/K$ modulo the ideal generated by a regular sequence preserves homological information: if $\mathcal{F}$ is a minimal free resolution of $T/K$, then $\mathcal{F} \otimes_T T/\langle A \rangle$ is a minimal free resolution of $R/I$ and~$\mathcal{F} \otimes_T T/\langle B \rangle$ is a minimal free resolution of $S/J$ (see~\cite[Proposition~1.1.5]{bruns}). Thus, these modules all have the same Betti numbers, and by Corollary~\ref{cor:cutbetti} the Betti numbers are given as above. $\square$

\vspace{-0.4cm}

\bibliography{wilmes}{}
\bibliographystyle{plain}

\end{document}